\documentclass[reqno]{amsart}
\usepackage{hyperref}

\begin{document}
\title{On the norming constants of the Sturm-Liouville problem}

\author{Tigran Harutyunyan, Avetik Pahlevanyan}

\address{Tigran Harutyunyan \newline
Faculty of Mathematics and Mechanics, Yerevan State University, 1 Alex Manoogian, 0025, Yerevan, Armenia}
\email{hartigr@yahoo.co.uk}
\address{Avetik Pahlevanyan \newline
Institute of Mathematics of National Academy of Sciences, 24/5 Baghramian ave., 0019, Yerevan, Armenia}
\email{apahlevanyan@instmath.sci.am}


\begin{abstract}
We derive new asymptotic formulae for the norming constants of Sturm-Liouville problem with summable potentials, which generalize and make more precise previously known formulae. Moreover, our formulae take into account the smooth dependence of norming constants on boundary conditions. We also find some new properties of the remainder terms of asymptotics.
\end{abstract}

\maketitle
\numberwithin{equation}{section}
\newtheorem{theorem}{Theorem}[section]
\newtheorem{lemma}[theorem]{Lemma}
\newtheorem{definition}[theorem]{Definition}
\newtheorem{remark}[theorem]{Remark}

\section{Introduction and Statement of the Results}
\label{sec1}

Let $L\left(q,\alpha ,\beta \right)$ denote the Sturm--Liouville problem
\begin{equation}\label{eq1.1}
- y'' + q\left(x\right)y = \mu y, \; x \in \left(0, \pi \right), \; \mu \in \mathbb{C},
\end{equation}
\begin{equation}\label{eq1.2}
y\left(0 \right)\cos \alpha + y'\left( 0 \right)\sin \alpha  = 0, \; \alpha \in \left(0, \pi\right],
\end{equation}
\begin{equation}\label{eq1.3}
y\left(\pi\right)\cos \beta + y'\left(\pi \right)\sin \beta  = 0, \; \beta \in \left[0, \pi\right),
\end{equation}
where $q$ is a real-valued, summable function on $\left[0, \pi\right]$ (we write $q \in L_\mathbb{R}^1\left[0, \pi \right]$). By $L\left(q, \alpha, \beta \right)$ we also denote the self-adjoint operator, generated by the problem \eqref{eq1.1}--\eqref{eq1.3} in Hilbert space $L^2\left[0, \pi\right]$ (see \cite{Naimark:1969,Levitan-Sargsyan:1970}). It is well-known that the spectra of $L\left(q, \alpha, \beta \right)$ is discrete and consists of real, simple eigenvalues (see \cite{Naimark:1969,Levitan-Sargsyan:1970,Marchenko:1977}), which we denote by $\mu_n\left(q, \alpha, \beta \right),$ $n=0,1,2\dots,$ emphasizing the dependence of $\mu_n$ on $q,$ $\alpha$ and $\beta.$ For $\mu_n$ the following asymptotic formula have been proven in \cite{Harutyunyan:2008}:
\begin{equation}\label{eq1.4}
\mu_n \left(q, \alpha, \beta \right) = \left[n+\delta_n \left(\alpha, \beta \right)\right]^2+\cfrac{1}{\pi} \displaystyle\int_{0}^{\pi}q\left(t\right)dt+r_n\left(q, \alpha, \beta \right),
\end{equation}
where $\delta_n$ is the solution of the equation (for $n \geq 2$)
\begin{multline}\label{eq1.5}
{\delta_n}(\alpha, \beta)=\cfrac{1}{\pi}\arccos \cfrac{\cos \alpha}{\sqrt {\left(n+{\delta_n}(\alpha, \beta)\right)^2\sin^2 \alpha+\cos^2 \alpha}}- \\
-\cfrac{1}{\pi}\arccos \cfrac{\cos \beta}{\sqrt {\left(n+{\delta_n}(\alpha, \beta)\right)^2\sin^2 \beta+\cos^2 \beta}},
\end{multline}
and $r_n \left(q, \alpha, \beta \right)=o\left(1\right),$ when $n \to \infty,$ uniformly in $\alpha, \beta \in \left[0, \pi \right]$ and $q$ from any bounded subset of $L_\mathbb{R}^1 \left[0, \pi \right]$ (we will write $q \in {BL}_\mathbb{R}^1 \left[0, \pi \right]$). It follows from \eqref{eq1.5} (for details see \cite{Harutyunyan:2008}), that
$$\delta_n \left(\alpha, \beta \right)=\cfrac{\cot \beta-\cot \alpha}{\pi n}+O\left( 1/n^2 \right), \; \alpha, \beta \in \left(0, \pi \right), \eqno (1.5\mbox{a})$$
$$\delta_n \left(\pi, \beta \right)=\cfrac{1}{2}+\cfrac{\cot \beta}{\pi \left(n+\frac{1}{2}\right)}+O\left( 1/n^2 \right)=\cfrac{1}{2}+O\left( 1/n \right), \; \beta \in \left(0, \pi \right), \eqno (1.5\mbox{b})$$
$$\delta_n \left(\alpha, 0 \right)=\cfrac{1}{2}-\cfrac{\cot \alpha}{\pi \left(n+\frac{1}{2}\right)}+O\left( 1/n^2 \right)=\cfrac{1}{2}+O\left( 1/n \right), \; \alpha \in \left(0, \pi \right), \eqno (1.5\mbox{c})$$
$$\delta_n \left(\pi, 0 \right)=1. \; \eqno (1.5\mbox{d})$$
Let $y=\varphi \left(x, \mu, \alpha, q\right)$ and $y=\psi \left(x, \mu, \beta, q\right)$ be the solutions of \eqref{eq1.1} with initial values
$$\varphi \left(0, \mu, \alpha, q \right)=\sin \alpha, \; \varphi' \left(0, \mu, \alpha, q \right)=-\cos \alpha,$$
$$\psi \left(\pi, \mu, \beta, q \right)=\sin \beta, \; \psi' \left(\pi, \mu, \beta, q \right)=-\cos \beta.$$
The eigenvalues $\mu_n$ of $L\left(q, \alpha, \beta \right)$ are the solutions of the equation
\begin{multline*}
\Phi \left(\mu \right)=\varphi \left(\pi, \mu, \alpha, q \right) \cos \beta + \varphi' \left(\pi, \mu, \alpha, q \right) \sin \beta=\\
=-\left[\psi \left(0, \mu, \beta, q \right) \sin \alpha + \psi' \left(0, \mu, \beta, q \right) \cos \alpha \right]=0.
\end{multline*}

It is easy to see that for arbitrary $n=0,1,2,\dots,$ \, $\varphi_n \left( x \right):=\varphi \left(x, \mu_n \left(q, \alpha, \beta \right), \alpha, q \right)$ and
$\psi_n \left( x \right):=\psi \left( x, \mu_n \left(q, \alpha, \beta \right), \beta, q \right)$ are eigenfunctions, corresponding to the eigenvalue $\mu_n \left( q, \alpha, \beta \right).$ The squares of the $L^2$-norm of these eigenfunctions:
$$a_n \left( q, \alpha, \beta \right):=\int_0^{\pi}{\left| \varphi_n \left( x \right) \right|^2 dx}, \;
  b_n \left( q, \alpha, \beta \right):=\int_0^{\pi}{\left| \psi_n \left( x \right) \right|^2 dx}$$
are called the norming constants.

The main results of this paper are the following theorems:

\begin{theorem}\label{thm1.1}
For norming constants $a_n$ and $b_n$ the following asymptotic formulae hold (when $n \to \infty$):
\begin{multline}\label{eq1.6}
a_n \left( q, \alpha, \beta \right) =\cfrac{\pi}{2} \left[ 1+ \cfrac {2 \, \ae_n \left( q, \alpha, \beta \right)}{\pi \left[n + \delta \left( \alpha, \beta \right)\right]} + r_n \right] \sin^2 \alpha + \\
+\cfrac{\pi}{2 \left[n + \delta_n(\alpha, \beta)\right]^2} \left[ 1+ \cfrac {2 \, \ae_n \left( q, \alpha, \beta \right)}{\pi \left[n + \delta \left( \alpha, \beta \right)\right]} + \tilde{r}_n \right]\cos^2 \alpha,
\end{multline}

\begin{multline*}
b_n \left( q, \alpha, \beta \right) =\cfrac{\pi}{2} \left[ 1+ \cfrac {2 \, \ae_n \left( q, \alpha, \beta \right)}{\pi \left[n + \delta \left( \alpha, \beta \right)\right]} + p_n \right] \sin^2 \beta + \\
+\cfrac{\pi}{2 \left[n + \delta_n(\alpha, \beta)\right]^2} \left[ 1+ \cfrac {2 \, \ae_n \left( q, \alpha, \beta \right)}{\pi \left[n + \delta \left( \alpha, \beta \right)\right]} + \tilde{p}_n \right]\cos^2 \beta,
\end{multline*}
where
\begin{equation}\label{eq1.7}
\ae_n = \ae_n \left( q, \alpha, \beta \right) = -\cfrac {1}{2} \displaystyle \int_{0}^{\pi} \left( \pi - t \right) q\left(t\right) \sin 2 \left[ n + \delta_n\left(\alpha, \beta\right)\right]t dt,
\end{equation}
$r_n = r_n \left( q, \alpha, \beta \right) = O \left( \cfrac {1} {n^2}\right)$ and $\tilde{r}_n=\tilde{r}_n \left( q, \alpha, \beta \right) = O \left( \cfrac {1} {n^2}\right)$ (the same estimate is true for $p_n$ and $\tilde{p}_n$), when $n \to \infty,$ uniformly in $\alpha, \beta \in \left[0, \pi\right]$ and $q \in {BL}^1_\mathbb{R}\left[0, \pi\right].$
\end{theorem}

\begin{theorem}\label{thm1.2}
For both $\alpha, \beta \in (0, \pi)$ and $\alpha = \pi, \; \beta = 0$ cases the function $k,$ defined as the series
\begin{equation*}
k(x) = \displaystyle \sum_{n=2}^{\infty} \cfrac {\ae_n}{n + \delta_n \left( \alpha, \beta \right)} \cos \left[n + \delta_n \left( \alpha, \beta \right)\right]x
\end{equation*}
is absolutely continuous function on arbitrary segment $\left[a,b\right] \subset \left(0, 2 \pi\right),$ i.e. $k \in AC \left(0, 2 \pi\right).$
\end{theorem}

The dependence of  norming constants on $\alpha$ and $\beta$ (as far as we know) hasn't been investigated before. The dependence of spectral data (by spectral data here we understand the set of eigenvalues and the set of norming constants) on $\alpha$ and $\beta$ has been usually studied (see \cite{Naimark:1969,Levitan-Sargsyan:1970,Marchenko:1977,Isaacson-Trubowitz:1983,Isaacson-McKean-Trubowitz:1984,
Dahlberg-Trubowitz:1984,Poschel-Trubowitz:1987,Yurko:2007}) in the following sense: the boundary conditions are separated into four cases:
\begin{itemize}
\item[1)]
$\sin \alpha \neq 0, \; \sin \beta \neq 0,$ \; i.e. $\alpha, \beta \in \left( 0, \pi \right);$
\item[2)]
$\sin \alpha = 0, \; \sin \beta \neq 0,$ \; i.e. $\alpha = \pi, \; \beta \in \left( 0, \pi \right);$
\item[3)]
$\sin \alpha \neq 0, \; \sin \beta = 0,$ \; i.e. $\alpha \in \left( 0, \pi \right), \; \beta = 0;$
\item[4)]
$\sin \alpha = 0, \; \sin \beta = 0,$ \; i.e. $\alpha = \pi, \; \beta = 0,$
\end{itemize}
and results are formulated separately for each case. For eigenvalues, formula \eqref{eq1.4} generalizes and unites four different formulae that were known before in four mentioned cases (see \cite{Harutyunyan:2008}).

So far, for norming constants the following is known.

In the case $ \sin \alpha \neq 0$ it is known that for smooth $q$
\begin{equation}\label{eq1.8}
\cfrac{a_n \left( q, \alpha, \beta \right)}{\sin^2 \alpha}=\cfrac{\pi}{2}+O \left(\cfrac{1}{n^2} \right).
\end{equation}
For absolutely continuous $q$ (we will write $q \in AC \left[ 0, \pi \right]$) the proof of \eqref{eq1.8} can be found in \cite{Levitan-Sargsyan:1970}. Let us note, that if $q \in AC\left[0, \pi\right],$ then $\ae_n = O \left(\cfrac{1}{n}\right),$ and it is easy to see, that in this case \eqref{eq1.6} takes the form \eqref{eq1.8}. In \cite{Zhikov:1967}, under the condition $q \left( x \right) = \cfrac{dF\left(x\right)}{dx}$ (almost everywhere), and $\sin\alpha \neq 0,$ where $F$ is a function of bounded variation (we will write $F \in BV\left[0, \pi\right]$), the author asserts that
\begin{equation}\label{eq1.9}
\cfrac{a_n \left( q, \alpha, \beta \right)}{\sin^2 \alpha} = \cfrac{\pi}{2}+\alpha_n,
\end{equation}
where the sequence $\left\{\alpha_n\right\}_{n = 0}^\infty$ is characterized by the condition that the function $f\left( x \right):=\displaystyle \sum\limits_{n = 0}^\infty \alpha _n \cos nx$ has a bounded variation on $\left[0, \pi \right],$ i.e. $f \in BV \left[ 0, \pi \right].$ Our result is similar to this, but there are some differences, in particular, we assert that $k \in AC \left(0, 2\pi\right).$

In \cite{Yurko:2007}, for $q \in L_\mathbb{R}^2 \left[0, \pi \right],$ it was proved that
\begin{equation}\label{eq1.10}
\cfrac{a_n \left( q, \alpha, \beta \right)}{\sin^2 \alpha} = \cfrac{\pi}{2} + \cfrac{\kappa_n}{n},
\end{equation}
where $\left\{\kappa_n \right\}_{n=0}^\infty \in l^2$ (i.e. $\displaystyle \sum\limits_{n = 0}^\infty \left| \kappa _n \right|^2 < \infty$), and $\kappa_n = \ae_n + O \left( \cfrac{1}{n} \right)$ (see \eqref{eq1.7}).

It is also important to note that norming constants $a_n\left(q, \alpha, \beta \right)$ are analytic functions on $\alpha$ and $\beta$. It easily follows from formulae \eqref{eq3.1}, \eqref{eq3.2} and \eqref{eq3.4} below and from the result in \cite{Harutyunyan:2008}, which states that $\lambda_n \left( q, \alpha, \beta \right)$ $\left(\lambda_n^2 \left( q, \alpha, \beta \right)=\mu_n \left( q, \alpha, \beta \right)\right)$ depend analytically on $\alpha$ and $\beta.$

In the case $\sin \alpha =0,$ $\sin \beta \neq 0$ it is known that for smooth $q$ (for $q \in AC \left[ 0, \pi \right]$ the proof of \eqref{eq1.11} can be found in \cite{Levitan-Sargsyan:1970})
\begin{equation}\label{eq1.11}
a_n \left( q, \pi, \beta \right) = \cfrac{\pi}{2 \left(n+1/2\right)^2} \left[1+O\left(\cfrac{1}{n^2}\right)\right].
\end{equation}
Since $\delta_{n} \left( \pi, \beta \right)=\cfrac{1}{2}+O \left( \cfrac{1}{n} \right)$ (see (1.5b)), then it is easy to see that \eqref{eq1.11} follows from \eqref{eq1.6}. Besides, we see that \eqref{eq1.6} smoothly turns into \eqref{eq1.11} when $\alpha \to \pi$ (for $q \in AC \left[0, \pi\right]$).

In the case $\sin \alpha =0$, $\sin \beta =0$ the following result can be found in \cite{Levitan-Sargsyan:1970} for $q \in AC \left[ 0, \pi \right]$:
\begin{equation*}
a_n \left(q, \pi, 0\right)=\cfrac{\pi}{2 n^{2}} \left[1+O\left(\cfrac{1}{n^2}\right)\right].
\end{equation*}

We think that it is more correct to write this result in the form (note that $\delta_n \left(\pi, 0 \right)=1$)
\begin{equation}\label{eq1.12}
a_n \left(q, \pi, 0\right)=\cfrac{\pi}{2 \left( n+1 \right)^2} \left[1+O\left(\cfrac{1}{n^2}\right)\right]
\end{equation}
to  keep the beginning of the enumeration of eigenvalues and norming constants starting from $0,$ but not from $1,$ as in \cite{Levitan-Sargsyan:1970}.

Our proofs of the theorems are based on the detailed study of the dependence of eigenfunctions $\varphi_{n}$ and $\psi_{n}$ on parameters $\alpha$ and $\beta.$ We will present it in the sections \ref{sec3} and \ref{sec4}. But first we need to prove some properties of the solutions of the equation \eqref{eq1.1}.

\section{Asymptotics of the solutions}
\label{sec2}

Let $q\in L_\mathbb{C}^1 \left[ 0, \pi \right],$ i.e. $q$ is a complex-valued, summable function on $\left[0, \pi\right],$ and let us denote by $y_i \left( x, \lambda \right),$ $i=1,2,3,4,$ the solutions of the equation
\begin{equation}\label{eq2.1}
-y''+q \left( x \right)y=\lambda^2 y,
\end{equation}
satisfying the initial conditions
\begin{equation}\label{eq2.2}
\begin{aligned}
y_1 \left(0, \lambda \right)=1, \; y_2 \left(0, \lambda \right)=0, \; y_3 \left(\pi, \lambda \right)=1, \; y_4 \left(\pi, \lambda \right)=0, \\
y'_1 \left(0, \lambda \right)=0, \; y'_2 \left(0, \lambda \right)=1, \; y'_3 \left(\pi, \lambda \right)=0, \; y'_4 \left(\pi, \lambda \right)=1.
\end{aligned}
\end{equation}

Let us recall that by a solution of \eqref{eq2.1} (which is the same as \eqref{eq1.1}) we understand the function $y,$ such that $y, \, y' \in AC \left[0, \pi \right]$ and which satisfies \eqref{eq2.1} almost everywhere (see \cite{Naimark:1969}).

The solutions $y_1$ and $y_2$ (as well as the second pair $y_3$ and $y_4$) form a fundamental system of  solutions of \eqref{eq1.1}, i.e. any solution $y$ of \eqref{eq1.1} can be represented in the form:
\begin{equation}\label{eq2.3}
y \left(x \right)=y \left(0 \right)y_1 \left(x, \lambda \right)+y' \left(0 \right)y_2\left(x, \lambda \right)=y \left(\pi \right)y_3 \left(x, \lambda \right)+y' \left(\pi \right)y_4 \left(x, \lambda \right).
\end{equation}

The existence and uniqueness of the solutions $y_i,$ $i=1,2,3,4$ (under the condition $q \in L_\mathbb{C}^1 \left[0, \pi\right]$) were investigated in \cite{Naimark:1969,Marchenko:1952,Chudov:1949,Atkinson:1964,Harutyunyan-Hovsepyan:2005}. The following lemma in some sense extends the results of the mentioned papers related to asymptotics (when $\left| \lambda \right| \to \infty$) of the solutions $y_i,$ $i=1,2,3,4.$

\begin{lemma}\label{lem2.1}
Let $q \in L_\mathbb{C}^1 \left[0, \pi\right].$ Then for the solutions $y_i,$ $i=1,2,3,4,$ the following representations hold (when $\left| \lambda \right| \geq 1$):
\begin{equation}\label{eq2.4}
y_1 \left(x, \lambda \right)=\cos \lambda x+\cfrac{1}{2\lambda} \, a \left(x, \lambda \right),
\end{equation}
\begin{equation}\label{eq2.5}
y_2 \left(x, \lambda \right)=\cfrac{\sin \lambda x}{\lambda}-\cfrac{1}{2 \lambda^2} \, b \left(x, \lambda \right),
\end{equation}
\begin{equation}\label{eq2.6}
y_3 \left(x, \lambda \right)=\cos \lambda \left(\pi -x \right)+\cfrac{1}{2\lambda} \, c \left(x, \lambda \right),
\end{equation}
\begin{equation}
y_4 \left(x, \lambda \right)=\cfrac{\sin \lambda \left(\pi -x \right)}{\lambda}-\cfrac{1}{2 \lambda^2} \, d\left(x, \lambda \right),
\end{equation}
where $a$, $b$, $c$, $d$ are twice differentiable with respect to $x$ and entire functions with respect to $\lambda,$ and have the form
\begin{equation}\label{eq2.8}
a\left(x, \lambda\right)= \sin \lambda x \displaystyle \int \limits^x_0 q\left(t\right) dt + \displaystyle \int \limits^x_0 q\left(t\right) \sin \lambda \left(x-2t\right) dt + R_1\left(x, \lambda, q\right),
\end{equation}
\begin{equation}\label{eq2.9}
b\left(x, \lambda\right)= \cos \lambda x \displaystyle \int \limits^x_0 q\left(t\right) dt - \displaystyle \int \limits^x_0 q\left(t\right) \cos \lambda \left(x-2t\right) dt + R_2\left(x, \lambda, q\right),
\end{equation}
\begin{equation}\label{eq2.10}
c\left(x, \lambda\right)= \sin \lambda \left(\pi -x\right) \displaystyle \int \limits^\pi_x q\left(t\right) dt + \displaystyle \int\limits^\pi_x q\left(t\right) \sin\lambda \left(2t-\pi-x\right) dt + R_3\left(x,\lambda, q\right),
\end{equation}
\begin{equation}\label{eq2.11}
d\left(x, \lambda\right)= \cos \lambda \left(\pi -x\right) \displaystyle \int \limits^\pi_x q\left(t\right) dt + \displaystyle \int \limits^\pi_x q\left(t\right) \cos \lambda \left(\pi+x-2t\right) dt + R_4\left(x,\lambda, q\right),
\end{equation}
and  $R_i,$ $i=1,2,3,4,$ satisfy the estimates (when $\left| \lambda \right| \geq 1$)
\begin{equation}\label{eq2.12}
R_1 \left(x, \lambda, q\right), \; R_2 \left(x, \lambda, q\right) = O \left( \cfrac{e^{\left| Im \lambda \right| x}}{\left| \lambda \right|} \right),
\end{equation}
\begin{equation}\label{eq2.13}
R_3 \left(x, \lambda, q\right), \; R_4 \left(x, \lambda, q\right) = O \left( \cfrac{e^{\left| Im \lambda \right| (\pi -x)}}{\left| \lambda \right|} \right),
\end{equation}
uniformly with respect to $q \in BL_{\mathbb{C}}^1 \left[0, \pi \right].$
\end{lemma}
\begin{proof} In \cite{Harutyunyan-Hovsepyan:2005} the authors have proved that $y_2 \left(x, \lambda \right)$ can be obtained as a sum of series
$$y_2 \left(x, \lambda, q \right)=\displaystyle \sum_{k=0}^{\infty}{S_k \left(x, \lambda, q \right)},$$
which converge to $y_2 \left(x, \lambda, q \right)$ uniformly on bounded subsets of the set $\left[0, \pi \right] \times \mathbb{C} \times L_\mathbb{C}^1 \left[0, \pi\right],$ and where $S_0 \left(x, \lambda, q \right)=\cfrac{\sin \lambda x}{\lambda},$
$$S_k \left(x, \lambda, q \right)=\displaystyle \int_0^x \cfrac{\sin \lambda \left(x-t \right)}{\lambda} \, q\left(t\right) S_{k-1} \left(t, \lambda, q \right)dt, \; k=1,2,\dots.$$
For $S_k$ we have  the estimate (when $\left| \lambda \right| \geq 1$):
\begin{equation}\label{eq2.14}
\left| S_k \left(x, \lambda, q \right) \right| \leq \cfrac{e^{\left| Im \lambda \right| x}}{{\left| \lambda \right|}^{k+1}} \; \cfrac{\sigma_0^k \left(x\right)}{k!}, \; k=0,1,2,\dots,
\end{equation}
where $\sigma_0 \left( x \right) \equiv \displaystyle \int_0^x \left| q \left(t \right) \right|dt$ (see \cite{Harutyunyan-Hovsepyan:2005}). To prove \eqref{eq2.5}, \eqref{eq2.9} and the estimate \eqref{eq2.12}, we write $S_1$ in the form
\begin{multline*}
S_1 \left(x, \lambda, q \right)=\displaystyle \int_0^x \cfrac{\sin \lambda \left(x-t \right)}{\lambda} \; \cfrac{\sin \lambda t}{\lambda} \, q \left(t \right)dt=\\
=\cfrac{1}{2 \lambda^2} \displaystyle \int_0^x \left[\cos \lambda \left(x-2t \right)-\cos \lambda x \right] q \left(t \right)dt=\\
=-\cfrac{\cos \lambda x}{2 \lambda^2} \displaystyle \int_0^x q \left(t \right)dt+\cfrac{1}{2 \lambda^2} \displaystyle \int_0^x{\cos \lambda \left(x-2t \right)} q \left(t \right)dt,
\end{multline*}
and note that
\begin{equation*}
S'_k \left(x, \lambda, q \right)=\displaystyle \int_0^x {\cos \lambda \left(x-t \right)} q \left(t \right) S_{k-1} \left(t, \lambda, q \right)dt, \; k=1,2,\dots.
\end{equation*}
This implies that $S'_k \in AC \left[0, \pi \right].$ By writing $y_2 \left(x, \lambda, q \right)=S_0+S_1+\displaystyle \sum_{k=2}^\infty S_k \left(x, \lambda, q\right),$ we obtain
\begin{equation*}
y_2 \left(x, \lambda, q \right)=\cfrac{\sin \lambda x}{\lambda}-\cfrac{1}{2 \lambda^2} \, b\left(x, \lambda \right),
\end{equation*}
where $-\cfrac{1}{2 \lambda^2} \, b \left(x, \lambda \right)=\displaystyle\sum_{k=1}^\infty S_k \left(x, \lambda, q\right)$ and therefore $b \left(x, \lambda \right)$ has the form \eqref{eq2.9}, where $R_2 \left(x, \lambda\right)= -2 \lambda^2 \displaystyle\sum_{k=2}^\infty S_k \left(x, \lambda\right).$ Now, from the estimate \eqref{eq2.14}, we obtain that
\begin{multline*}
\sum_{k=2}^\infty \left|S_k \left(x, \lambda, q \right) \right| \leq \sum_{k=2}^\infty \cfrac{e^{\left| Im \lambda \right|x}}{\left| \lambda \right|^{k+1}} \; \cfrac{\sigma_0^k \left(x \right)}{k!}=\cfrac{e^{\left| Im \lambda \right|x} \, \sigma_0^2 \left(x \right)}{\left| \lambda \right|^3} \sum_{k=2}^{\infty} \cfrac{\sigma_0^{k-2} \left(x \right)}{\left| \lambda \right|^{k-2} k!}< \\
<\cfrac{e^{\left| Im \lambda \right|x} \, \sigma_0^2 \left(x \right)}{\left| \lambda \right|^3} \sum_{k=2}^{\infty} \cfrac{\sigma_0^{k-2} \left(x \right)}{\left| \lambda \right|^{k-2} \left(k-2 \right)!}=\cfrac{e^{\left| Im \lambda \right|x} \, \sigma_0^2 \left(x\right)}{\left|\lambda \right|^3} \sum_{n=0}^{\infty} \cfrac{\sigma_0^n \left(x \right)}{\left| \lambda \right|^n n!}= \\
=\cfrac{e^{\left| Im \lambda \right|x} \, \sigma_0^2 \left(x \right)}{\left| \lambda \right|^3} \;  e^{\tfrac{\sigma_0 \left(x \right)}{\left| \lambda \right|}}=\cfrac{\sigma_0^2 \left(x \right)}{\left| \lambda \right|^3} \; e^{\left| Im \lambda \right|x+\tfrac{\sigma_0 \left(x \right)}{\left|\lambda \right|}}.
\end{multline*}
This implies \eqref{eq2.12} for $\left| \lambda \right| \geq 1$. Since $y'_2 \in AC \left[ 0, \pi \right]$ and $S_0 \left( x, \lambda \right)=\cfrac{\sin \lambda x}{\lambda},$ then we obtain that $-\cfrac{1}{2 \lambda^2} \,
b \left( x, \lambda \right)= y_2 - S_0$ is also a twice differentiable function (more precisely $b' \in AC \left[ 0, \pi \right]$). Assertions for $y_1,$ $y_3,$ $y_4$ can be proven similarly.
\end{proof}

\section{The proof of the Theorem \ref{thm1.1}}
\label{sec3}

According to \eqref{eq2.3}, the solution $\varphi \left(x, \mu, \alpha, q \right),$ which we will denote by $\varphi \left(x, \lambda^2, \alpha \right)$ for brevity, has the form
\begin{equation}\label{eq3.1}
\varphi \left(x, \lambda^2, \alpha \right)=y_1 \left(x, \lambda \right) \sin \alpha - y_2 \left(x, \lambda \right) \cos \alpha,
\end{equation}
and according to \eqref{eq2.4} and \eqref{eq2.5} we arrive at:
\begin{equation*}
\varphi \left(x, \lambda^2, \alpha \right)=\left[\cos \lambda x+\cfrac{1}{2\lambda} \, a \left(x, \lambda \right) \right] \sin \alpha - \left[\cfrac{\sin \lambda x}{\lambda}-\cfrac{1}{2 \lambda^2} \, b \left(x, \lambda \right) \right] \cos \alpha.
\end{equation*}
Taking the squares of both sides of the last equality, we obtain:
\begin{multline}\label{eq3.2}
\varphi^2 \left(x, \lambda^2, \alpha \right)=\cos^2 \lambda x \sin^2 \alpha + \cfrac{1}{\lambda} \left[a\left(x,\lambda \right) \cos \lambda x + \cfrac{a^2 \left(x, \lambda \right)}{4\lambda}\right] \sin^2 \alpha - \\
-\cfrac{2}{\lambda} \left[\cos \lambda x \sin \lambda x-\cfrac{b\left(x, \lambda \right) \cos \lambda x} {2\lambda}+\cfrac{a\left(x, \lambda \right) \sin \lambda x}{2\lambda}-\cfrac{a\left(x, \lambda \right) b\left(x, \lambda \right)}{4 \lambda^2} \right] \times \\
\times \sin \alpha \cos \alpha +\cfrac{\sin^2 \lambda x}{\lambda^2} \, \cos^2 \alpha +\left[\cfrac{b^2 \left(x, \lambda \right)}{4 \lambda^4}-\cfrac{b\left(x, \lambda \right) \sin \lambda x}{\lambda^3}\right] \cos^2 \alpha. \end{multline}
Recalling the formulae $\cos^2 \lambda x = \cfrac{1}{2}\left(1+\cos 2 \lambda x \right)$ and $\sin^2 \lambda x = \cfrac{1}{2} \left(1-\cos 2 \lambda x \right),$ from \eqref{eq3.2}, we obtain:
\begin{multline}\label{eq3.3}
\displaystyle \int_0^{\pi} \varphi^2 \left(x, \lambda^2, \alpha \right)dx = \cfrac{\pi}{2} \sin^2 \alpha +\cfrac{\sin 2\lambda \pi}{4\lambda} \, \sin^2 \alpha + \\
+\cfrac{1}{\lambda}\left(\int_0^{\pi} a\left(x, \lambda \right) \cos \lambda xdx + \cfrac{1}{4\lambda} \int_0^{\pi} a^2 \left(x, \lambda \right)dx \right) \sin^2 \alpha-\cfrac{\sin^2 \lambda \pi}{\lambda^2} \sin \alpha \cos \alpha + \\
+\cfrac{1}{\lambda^2} \left(\displaystyle\int_0^{\pi} b\left(x, \lambda \right) \cos \lambda xdx - \displaystyle\int_0^{\pi} a\left(x, \lambda \right) \sin \lambda xdx \right) \sin \alpha \cos \alpha + \\
+\cfrac{\sin \alpha \cos \alpha}{2\lambda^3} \int_0^{\pi} a\left(x, \lambda \right) b\left(x, \lambda \right) dx + \cfrac{\pi}{2 \lambda^2} \, \cos^2 \alpha - \cfrac{\sin 2\lambda \pi}{4 \lambda^3} \, \cos^2 \alpha - \\
-\cfrac{1}{\lambda^3}\left(\displaystyle\int_0^{\pi} b\left(x, \lambda \right) \sin \lambda xdx - \cfrac{1}{4\lambda}\displaystyle\int_0^{\pi} b^2 \left(x, \lambda \right)dx \right) \cos^2 \alpha.
\end{multline}
We are going to receive the asymptotic formula \eqref{eq1.6} by the substitution $\lambda=\lambda_n \left(q, \alpha, \beta \right)=\sqrt{\mu_n \left(q, \alpha, \beta \right)}$ in \eqref{eq3.3}. To this aim, we estimate each term of the right-hand side of \eqref{eq3.3} for $\lambda ={\lambda}_n$. It can be easily deduced from \eqref{eq1.4} that for $\lambda_n = \sqrt{\mu_n}$ we have the following asymptotic formula:
\begin{equation}\label{eq3.4}
\lambda_{n} \left( q, \alpha, \beta \right) = n+ \delta_{n} \left( \alpha, \beta \right)+ \cfrac{\left[ q \right]}{2{\left( n+\delta_{n} \left( \alpha, \beta \right) \right)}}+l_{n},
\end{equation}
where $ \left[ q \right] := \cfrac{1}{\pi} \displaystyle\int_0^{\pi} q \left(t\right)dt,$ $l_n=l_n \left(q, \alpha, \beta \right)=o \left(\cfrac{1}{n} \right)$ uniformly with respect to $\alpha, \beta \in \left[0, \pi \right]$ and $q \in {BL}_\mathbb{R}^1 \left[0,\pi \right]$ (see \cite{Harutyunyan:2016}).

It follows from (1.5a)--(1.5d) that $\sin 2\pi \delta_n \left(\alpha, \beta \right)=O\left(\cfrac{1}{n}\right)$ and
\begin{equation}\label{eq3.5}
\sin 2\pi \lambda_n \left(q, \alpha, \beta\right) = O\left(\cfrac{1}{n} \right), \; \cos 2\pi \lambda_n \left( q, \alpha, \beta \right) = 1 - O\left(\cfrac{1}{n^2} \right),
\end{equation}
for all $\left(\alpha, \beta\right) \in \left(0, \pi \right] \times \left[0, \pi \right).$

Thus, the second term
\begin{equation}\label{eq3.6}
\cfrac{\sin 2\pi \lambda_n}{\lambda_n} \, \sin^2 \alpha = O\left( \cfrac{1}{n^2} \right) \sin^2 \alpha.
\end{equation}
Important is the third term: $\cfrac{1}{\lambda_n} \displaystyle \int_0^{\pi} a\left(x, \lambda_n \right) \cos  \lambda_n x dx.$ According to \eqref{eq2.8} and \eqref{eq2.12} we have
\begin{equation}\label{eq3.7}
a \left(x, \lambda_n\right) = A \left(x, \lambda_n\right) + O \left( \cfrac {1}{\lambda_n}\right),
\end{equation}
where
\begin{equation}\label{eq3.8}
A \left(x, \lambda_n\right) = \displaystyle \int_0^x q\left(t\right)dt \sin \lambda_n x + \int_0^x q\left(t\right)\sin \lambda_n \left(x - 2t\right) dt.
\end{equation}
After multiplying both sides by $\cos \lambda_n x,$ integrating over $\left[0,\pi\right]$ and changing the order of integration we get
\begin{multline}\label{eq3.9}
\displaystyle \int_0^{\pi} A\left(x, \lambda_n\right) \cos \lambda_n x dx =
\cfrac{\sin^2 \lambda_n \pi}{\lambda_n} \displaystyle \int_0^{\pi} q\left(t\right) \cos^2 \lambda_n t dt - \\
- \cfrac{\sin 2 \lambda_n \pi}{4 \lambda_n} \displaystyle \int_0^{\pi} q\left(t\right) \sin 2 \lambda_n t dt - \cfrac{1}{2} \displaystyle \int_0^{\pi} \left(\pi - t\right) q\left(t\right) \sin 2 \lambda_n tdt.
\end{multline}
Taking into account the formulae \eqref{eq3.5} and denoting (see \eqref{eq1.7})
\begin{equation*}
\tilde{\ae}_n \equiv \tilde{\ae}_n \left( q, \alpha, \beta \right) := -\cfrac {1}{2} \displaystyle \int_0^{\pi} \left( \pi - t \right) q\left(t\right) \sin 2 \lambda_n tdt,
\end{equation*}
we can rewrite \eqref{eq3.9} in the form
\begin{equation}\label{eq3.10}
\displaystyle \int_{0}^{\pi} A(x, \lambda_n) \cos \lambda_n x dx = \tilde {\ae}_n + O \left( \cfrac{1}{n} \right).
\end{equation}
Since $\sin 2 \lambda_n t = \sin 2 \left(n+\delta_n + O \left( \cfrac{1}{n} \right)\right)t = \sin 2 \left(n+\delta_n \right)t + O \left( \cfrac{1}{n} \right)$ holds uniformly with respect to $t \in \left[0, \pi \right],$ then $\tilde{\ae}_n = \ae_n + O \left( \cfrac{1}{n} \right),$ and therefore the third term of \eqref{eq3.3} has the form
\begin{equation*}
\cfrac{1}{\lambda_n} \displaystyle \int_{0}^{\pi} a(x, \lambda_n) \cos \lambda_n x dx=\cfrac{\ae_n}{n + \delta_n (\alpha, \beta)} + O \left( \cfrac{1}{n^2} \right).
\end{equation*}
Now, let us focus on the remained terms of the equality \eqref{eq3.3} for $\lambda=\lambda_n.$ The terms from the fourth to the eighth  have the coefficient $\cfrac{1}{\lambda_n^{\gamma}},$ where $\gamma \geq 2,$ and therefore they have the order $O \left( \cfrac{1}{n^2} \right).$ Concerning the last four terms of \eqref{eq3.3}, we observe that both $\cfrac{\sin 2 \pi \lambda_n}{4 \lambda_n^3} \cos^2 \alpha$ and $\cfrac{1}{\lambda_n^4} \displaystyle \int_{0}^{\pi} b^2 \left(x, \lambda_n\right) dx \cos^2 \alpha$ have the same order $O \left( \cfrac{1}{\lambda_n^4} \right) \cos^2 \alpha.$ An important term is $\cfrac{1}{\lambda_n^3} \displaystyle \int_{0}^{\pi} b\left(x, \lambda_n\right) \sin \lambda_n x dx.$ According to \eqref{eq2.9} and \eqref{eq2.12} we can write $b \left(x, \lambda_n\right)$ in the form
\begin{equation}\label{eq3.11}
b \left(x, \lambda_n \right) = B \left(x, \lambda_n \right) + O \left( \cfrac{1}{\lambda_n} \right),
\end{equation}
where
\begin{equation}\label{eq3.12}
B\left(x, \lambda_n\right) =  \displaystyle \int_{0}^{x} q\left(t\right) dt \cos \lambda_n x - \displaystyle \int_{0}^{x} q\left(t\right) \cos \lambda_n \left(x-2t\right) dt.
\end{equation}
A simple computation yields:
\begin{equation*}
B\left(x, \lambda_n \right) \sin \lambda_n x =  \displaystyle \int_{0}^{x} q\left(t\right) dt \sin 2 \lambda_n x - \displaystyle \int_{0}^{x} q\left(t\right) \sin 2 \lambda_n t dt - A\left(x, \lambda_n \right) \cos \lambda_n x.
\end{equation*}
After integrating the latter equality from $0$ to $\pi,$ changing the order of integration and taking into consideration \eqref{eq3.9} we get:
\begin{multline*}
\displaystyle \int_{0}^{\pi} B\left(x,\lambda_n\right) \sin \lambda_n x dx = - \cfrac{\cos^2 \lambda_n \pi}{\lambda_n} \displaystyle \int_{0}^{\pi} q\left(t\right) \sin^2 \lambda_n t dt + \\
+\cfrac{\sin 2 \lambda_n \pi}{4 \lambda_n} \displaystyle \int_{0}^{\pi} q\left(t\right) \sin 2 \lambda_n t dt
 - \cfrac{1}{2} \displaystyle \int_{0}^{\pi} \left(\pi - t\right) q\left(t\right) \sin 2 \lambda_n t dt = \\
= O \left( \cfrac{1}{\lambda_n} \right) + O \left( \cfrac{1}{\lambda_n^2} \right) + \tilde{\ae}_n=\ae_n+O\left(\cfrac{1}{n}\right),
\end{multline*}
and therefore the eleventh term of the equality \eqref{eq3.3} for $\lambda=\lambda_n$ has the form
\begin{equation*}
\cfrac{1}{\lambda_n^3}\displaystyle\int_0^{\pi} b\left(x, \lambda_n \right) \sin \lambda_n xdx=\cfrac{1}{\lambda_n^2} \left(\cfrac{\ae_n}{n+\delta_n \left(\alpha,\beta\right)}+O\left(\cfrac{1}{n^2}\right)\right).
\end{equation*}
Let us remark that from \eqref{eq3.4} we have $\cfrac{1}{\lambda_n} - \cfrac{1}{n+\delta_n} = O \left( \cfrac{1}{n^3} \right).$ Thus,
\begin{multline}\label{eq3.13}
a_n \left(q, \alpha, \beta \right) = \cfrac{\pi}{2} \left[ 1+ \cfrac{2 \, \ae_n}{\pi \left[n+\delta_n \left(\alpha,\beta\right)\right]} + O \left( \cfrac{1}{n^2} \right) \right] \sin^2 \alpha + O \left( \cfrac{1}{n^2} \right) \sin \alpha \cos \alpha+\\
+\cfrac{\pi}{2 \left[n+\delta_n \left(\alpha,\beta\right)\right]^2} \left[ 1+ \cfrac{2 \, \ae_n}{\pi \left[n+\delta_n \left(\alpha,\beta\right)\right]} + O \left( \cfrac{1}{n^2} \right) \right] \cos^2 \alpha.
\end{multline}

If $\sin \alpha \neq 0,$ then $O \left( \cfrac{1}{n^2} \right) \sin \alpha \cos \alpha$ can be included into the term $O \left( \cfrac{1}{n^2} \right) \sin^2 \alpha,$ and if $\sin \alpha = 0,$ then these terms are absent. Finally, we can write \eqref{eq3.13} in the form \eqref{eq1.6}. For $b_n$ everything can be done similarly. Theorem \ref{thm1.1} is proved.

\section{The proof of the Theorem \ref{thm1.2}}
\label{sec4}

In the sequel the following notations will be used:
\begin{equation}\label{eq4.1}
\tilde{q} \left(t \right):= \left(\pi-t\right)q\left(t\right) \;\; \mbox{and} \;\; \sigma\left(x\right) := \displaystyle \int_{0}^{x} \tilde{q} \left(t \right) dt = \displaystyle \int_{0}^{x} \left(\pi - t \right) q \left(t \right) dt.
\end{equation}
Now, we have
\begin{multline}\label{eq4.2}
\cfrac{\ae_n}{n+\delta_n\left(\alpha, \beta\right)} = -\cfrac{1}{2 \left[n+\delta_n \left(\alpha, \beta\right)\right]} \displaystyle \int_{0}^{\pi} \tilde {q} \left(t\right) \sin 2 \left( n+\delta_n \right)t
dt = \\
=-\cfrac{1}{2 \left(n+\delta_n\right)} \displaystyle \int_{0}^{\pi} \sin 2 \left(n+\delta_n\right) t d \sigma \left(t\right) = -\cfrac{\sigma \left(\pi\right) \sin 2 \pi \delta_n}{2 \left(n+\delta_n \right)} + \\
+\displaystyle \int_{0}^{\pi} \sigma\left(t\right) \cos 2 \left(n+\delta_n \right)t dt.
\end{multline}
It was observed in \eqref{eq3.5} that $\sin 2 \pi \delta_n = O \left( \cfrac{1}{n} \right)$. If we denote by $\tilde{\sigma} \left(x\right) := \sigma\left(\cfrac{x}{2}\right)$ and $c_n := \cfrac {\sin 2 \pi \delta_n}{2 \left(n+\delta_n\right)}=O\left(\cfrac{1}{n^2}\right),$ then we can rewrite $k\left(x \right)$ in the form
\begin{equation}\label{eq4.3}
k\left(x\right) = k_1 \left(x\right) + k_2 \left(x\right),
\end{equation}
where
\begin{equation}\label{eq4.4}
k_1 \left(x\right) = -\sigma\left(\pi\right) \displaystyle \sum_{n=2}^{\infty} c_n \cos \left[n+\delta_n \left(\alpha, \beta \right) \right]x,
\end{equation}
\begin{equation}\label{eq4.5}
k_2 \left(x\right) = \displaystyle \sum_{n=2}^{\infty} \displaystyle \int_{0}^{2 \pi} \tilde{\sigma}\left(t\right) \cos \left[ n+\delta_n \left(\alpha, \beta \right) \right]t dt \cos \left[n+\delta_n \left(\alpha, \beta \right) \right]x.
\end{equation}
Since $c_n=O\left(\cfrac{1}{n^2}\right),$ then the series in \eqref{eq4.4} converges absolutely and uniformly on $\left[0,2\pi \right],$ and $k_1 \in AC \left[0,2\pi \right].$

Next, we consider two cases:

{\bf Case I:} If $\alpha, \beta \in (0, \pi),$ then by (1.5a) we have
\begin{equation*}
\delta_n(\alpha, \beta)=\cfrac{\cot \beta - \cot \alpha}{\pi n} + O \left( \cfrac{1}{n^2} \right)= \cfrac{d}{n} + O \left( \cfrac{1}{n^2}\right)=O \left( \cfrac{1}{n}\right),
\end{equation*}
where $d=\cfrac{\cot \beta - \cot \alpha}{\pi}.$
Recalling the Maclaurin expansions of the functions $\sin x$ and $\cos x$ around the point $x=0,$ we obtain
\begin{equation}\label{eq4.6}
\cos \left[n+\delta_n \left(\alpha, \beta\right)\right]x= \cos nx - d \cdot x \, \cfrac{\sin nx}{n}+ e_n \left(x\right),
\end{equation}
where $e_n \left(x\right),$ as all the other entries of \eqref{eq4.6}, is a smooth function $\left(e_n \in C^{\infty}\right)$ and
\begin{equation}\label{eq4.7}
e_n \left(x\right)=O\left(\cfrac{1}{n^2}\right)
\end{equation}
uniformly on $x \in \left[0,2\pi \right].$ Therefore $k_2$ can be written in the form
\begin{equation*}
k_2 \left(x\right) = l_1 \left(x\right)+l_2 \left(x\right)+l_3 \left(x\right),
\end{equation*}
where
\begin{multline}\label{eq4.8}
l_1 \left(x\right)=-d \cdot x \displaystyle \sum_{n=2}^{\infty} \cfrac{1}{n} \int_{0}^{2 \pi} \tilde{\sigma}\left(t\right) \cos nt dt \sin nx - \\
-d \cdot \displaystyle \sum_{n=2}^{\infty} \cfrac{1}{n} \int_{0}^{2 \pi} t \tilde{\sigma}\left(t\right) \sin nt dt \cos nx + d^2 \cdot x \displaystyle \sum_{n=2}^{\infty} \cfrac{1}{n^2} \int_{0}^{2 \pi} t \tilde{\sigma}\left(t\right) \sin nt dt \sin nx + \\
+\displaystyle \sum_{n=2}^{\infty} \int_{0}^{2 \pi} e_n \left(t \right) \tilde{\sigma}\left(t\right) dt \cos nx -d \cdot x \displaystyle \sum_{n=2}^{\infty} \cfrac{1}{n} \int_{0}^{2 \pi} e_n \left(t \right) \tilde{\sigma}\left(t\right) dt \sin nx,
\end{multline}
\begin{multline}\label{eq4.9}
l_2 \left(x\right)=\displaystyle \sum_{n=2}^{\infty} e_n \left(x \right) \int_{0}^{2 \pi} \tilde{\sigma}\left(t\right) \cos nt dt - \\
-d \cdot \displaystyle \sum_{n=2}^{\infty} \cfrac{e_n \left(x \right)}{n} \int_{0}^{2 \pi} t \tilde{\sigma}\left(t\right) \sin nt dt + \displaystyle \sum_{n=2}^{\infty} e_n \left(x \right) \int_{0}^{2 \pi} e_n \left(t \right) \tilde{\sigma}\left(t\right) \sin nt dt,
\end{multline}
\begin{equation}\label{eq4.10}
l_3 \left(x\right)=\displaystyle \sum_{n=2}^{\infty}\displaystyle \int_{0}^{2 \pi} \tilde{\sigma}\left(t\right) \cos nt dt \cos nx.
\end{equation}

Since $\tilde \sigma \in AC [0, 2 \pi],$ then Fourier coefficients are
\begin{equation}\label{eq4.11}
\displaystyle \int_{0}^{2 \pi} \tilde{\sigma}\left(t\right) \cos nt dt = O \left( \cfrac{1}{n} \right), \;
\displaystyle \int_{0}^{2 \pi} t\tilde{\sigma}\left(t\right) \sin nt dt = O \left( \cfrac{1}{n} \right).
\end{equation}
Also we note that
\begin{equation}\label{eq4.12}
\displaystyle \int_{0}^{2 \pi} e_n \left(t \right) \tilde{\sigma}\left(t\right) dt=O\left(\cfrac{1}{n^2}\right).
\end{equation}
Therefore the trigonometric series in \eqref{eq4.8} converges absolutely and uniformly on $\left[0,2\pi \right],$ and $l_1 \in AC \left(0,2\pi \right).$

It follows from \eqref{eq4.7}, \eqref{eq4.11} and \eqref{eq4.12} that the terms of the series in \eqref{eq4.9} have the order $O\left(\cfrac{1}{n^3}\right),$ and therefore $l_2 \in AC \left[0, 2\pi \right].$

About $l_3 \left(x\right)$ we can say the following:
Since $\tilde{\sigma} \in AC \left[0, 2 \pi\right],$ then the Fourier series of $\tilde{\sigma}$
\begin{equation*}
\tilde{\sigma}\left(x\right) =
  \cfrac{a_0 \left(\tilde{\sigma}\right)}{2} + \displaystyle \sum_{n=1}^{\infty}\left( a_n \left(\tilde{\sigma}\right) \cos nx + b_n \left(\tilde{\sigma}\right) \sin x\right),
\end{equation*}
where $a_n \left(\tilde{\sigma}\right)=\cfrac{1}{\pi} \displaystyle \int_{0}^{2\pi} \tilde{\sigma} \left(t\right) \cos ntdt,$ $b_n \left(\tilde{\sigma}\right)=\cfrac{1}{\pi} \displaystyle \int_{0}^{2\pi} \tilde{\sigma}\left(t\right) \sin nt dt,$
converges to $\tilde{\sigma}\left(x\right)$ in every point of $\left[0, 2 \pi\right]$ and this series is a function from $AC\left[0, 2 \pi\right].$
The same is true for $\sigma^{\star} (x) = \tilde{\sigma}\left(2 \pi - x\right):$
\begin{equation*}
\tilde{\sigma}\left(2 \pi - x\right) = \cfrac{a_0 \left(\sigma^{\star}\right)}{2} + \displaystyle \sum_{n=1}^{\infty}\left( a_n \left(\sigma^{\star}\right) \cos nx + b_n \left(\sigma^{\star}\right) \sin x \right).
\end{equation*}
But it is easy to see, that $a_n \left(\sigma^{\star}\right) = a_n \left(\tilde{\sigma}\right)$ and $b_n \left(\sigma^{\star}\right) = -b_n \left(\tilde{\sigma}\right).$ So
\begin{equation*}
\cfrac{1}{2} \left( \tilde{\sigma}\left(x\right) + \tilde{\sigma}\left(2\pi - x\right) \right) = \cfrac{a_0 \left(\tilde{\sigma}\right)}{2} + \displaystyle \sum_{n=1}^{\infty} a_n \left(\tilde{\sigma}\right) \cos nx,
\end{equation*}
i.e. this is "the even part" of Fourier series of $\tilde{\sigma} \left(x\right),$ and is absolutely continuous on $\left[0, 2\pi\right].$
Thus, for the case $\alpha, \beta \in \left(0, \pi\right)$ Theorem \ref{thm1.2} is proved.

{\bf Case II:} If $\alpha = \pi, \beta = 0,$ then $\delta_n \left(\pi, 0 \right) = 1,$ and the function $k_2 \left(\cdot\right)$ takes the form $k_2 \left(x\right) = \displaystyle \sum_{n=3}^{\infty} \int_{0}^{2 \pi} \tilde{\sigma} \left(t\right) \cos nt dt \cos nx$ and again it is "the even part" of Fourier series (without the zeroth, the first and the second terms) of an absolutely continuous function. Theorem \ref{thm1.2} is proved.

\subsection*{Acknowledgment}
T.N.~Harutyunyan was supported by State Committee of Science MES RA, in frame of the research project No. 15T--1A392.

\end{document}